\documentclass[a4paper,11pt]{article}
\usepackage[english]{babel}
\usepackage{graphicx}
\usepackage{bbding}
\usepackage[latin1]{inputenc}
\usepackage{fancyhdr}
\usepackage{amsmath,amssymb,color,epsfig}
\usepackage{mathrsfs}
\usepackage{eufrak}
\usepackage{dsfont}
\usepackage{upgreek}
\usepackage{fancyhdr}
\usepackage[dvipsnames]{xcolor}
\newtheorem{theorem}{Theorem}[section]
\newtheorem{proposition}[theorem]{Proposition}

\newtheorem{remark}[theorem]{Remark}

\def\1{\mathds{1}}

\title{Double Inequalities for Complete Monotonicity Degrees of Remainders of
Asymptotic Expansions of the Gamma and Digamma Functions}
\author{Mohamed Bouali }                 

\date{}
\begin{document}

\maketitle
\begin{abstract} Motivated by several conjectures posed in the paper " Completely monotonic degrees for a difference between the logarithmic
and psi functions",we confirm in this work some conjectures on completely monotonic degrees of remainders of the asymptotic expansion of the logarithm of the gamma function and the digamma function and we give two bounded for this degrees.
\end{abstract}
\section{Introduction}
Completely monotonic functions have attracted the attention of many authors. Mathematicians
have proved many interesting results on this topic. For example, Koumandos \cite{kou1} obtained upper and
lower polynomial bounds for the function $x/(e^x -1)$, $x > 0$, with coefficients of the Bernoulli numbers
$B_k$. This enabled him to give simpler proofs of some results of H. Alzer and F. Qi et al., concerning
complete monotonicity of certain functions involving the functions $\Gamma(x)$, $\psi(x)$ and the polygamma
functions $\psi^{(n)}(x)$, $n = 1, 2, ....$, \cite{clar}.

A function $f$ is said to be completely monotonic on an interval $I$ if $f$ has derivatives of all
orders on $I$ which alternate successively in sign, that is, $(-1)^nf^{(n)}(x)\geq 0$ for all $x\in I$ and all $n\in\Bbb N$. See for example [\cite{mit}, Chap VIII], [\cite{sch}, Chap I], and [\cite{wid}, Chap IV].

A notion of completely monotonic degree was invented first in reference \cite{guo1} and reviewed
in the recent paper \cite{Qi5}. It can be used to measure and differentiate complete
monotonicity more accurately, and it is also introduced in \cite{guo1,kou2,kou3,kou4,kou5,kou6, Qi1, Qi2, Qi3, Qi4} and closely related references.

Let $f(x)$ be a completely monotonic function on $(0,+\infty)$ and denote $f (+\infty) = \lim_{x\to +\infty}f(x)\geq 0$. When the function
$x^r[f(x)-f(+\infty)]$ is completely monotonic on $(0,+\infty)$ if and only if $0 \leq r\leq \alpha$, the number $\alpha$, denoted by ${\rm deg}^x_{cm}[f]$, is called the completely monotonic degree of $f(x)$ with respect to $x\in (0,+\infty)$. For more studies on complete monotonicity, the reader is also referred to \cite{guo1, kou2,kou3, kou6, Qi2, Qi3}.

For $x > 0$, the classical gamma function $\displaystyle\Gamma(x)=\int_0^\infty t^{x-1}e^{-t}dt$ first introduced by L. Euler, is one of the most important functions in mathematical analysis. It often appears in asymptotic series, hypergeometric series, Riemann zeta function, number theory, and so on.

In [\cite{alz2}, Theorem 8], [\cite{kou7}, Theorem 2], and \cite{xu}, the functions $$R_n(x)=(-1)^n[\log\Gamma(x)-(x-\frac12)\log(x)+x+\frac12\log(2\pi)-\sum_{k=1}^n\frac{B_{2k}}{2k(2k-1)}\frac1{x^{2k-1}},$$
for $n\geq 0$ were proved to be completely monotonic on $(0,+\infty)$, where an empty sum
is understood to be $0$ and the Bernoulli numbers $B_n$ are define by the following series \cite{Qi6, Qi8, Qi9} by
$$\frac z{e^z-1}=1-\frac z2+\sum_{n=1}^\infty\frac{b_{2n}}{(2n)!}z^{2n},\qquad|z|<2\pi.$$
Which implies that the functions $(-1)^mR_n^{(m)}$ for $m, n \geq 0$ are completely monotonic on $(0,+\infty)$. By the way, we call the function $(-1)^nR_n(x)$ for $n \geq 0$ the remainders of asymptotic
formula of $\log \Gamma(x)$. See [\cite{abr}, p. 257, 6.1.40] and [\cite{olv}, p. 140, 5.11.1]. The completely
monotonic degree of the function $R_n(x)$ for $n\geq 0$ with respect to $x$, $(0,+\infty)$ was
proved in [\cite{kou5}, Theorem 2.1] to be at least $n$.

Stimulated by the above results and related ones, Qi conjectured in [31] that:
the completely monotonic degrees of $R_n(x)$ for $n \geq 0$ with respect to $x$,
$(0,+\infty)$ satisfy
\begin{equation}\label{R1}{\rm deg}^x_{cm}(R_0)=0,\qquad {\rm deg}^x_{cm}(R_1)=1,\end{equation}
and
\begin{equation}\label{R2}{\rm deg}^x_{cm}(R_n)=2(n-1),\quad n\geq 2.\end{equation}
The completely monotonic degrees of $-R'_n(x)$ for $n \geq 0$ with respect to $x$,
\begin{equation}\label{R3}{\rm deg}^x_{cm}(-R'_0)=1,\qquad {\rm deg}^x_{cm}(-R'_1)=2.\end{equation}
and
\begin{equation}\label{R4}{\rm deg}^x_{cm}(-R'_n)=2n-1,\quad n\geq 2.\end{equation}
The completely monotonic degrees of $(-1)^mR^{(m)}_n(x)$ for $m\geq 2$ and $n \geq 0$ with respect to $x$,
\begin{equation}\label{R6}{\rm deg}^x_{cm}((-1)^m R_0^{(m)})=m-1,\quad {\rm deg}^x_{cm}((-1)^m R_1^{(m)})=m,\end{equation}
and
\begin{equation}\label{R5}{\rm deg}^x_{cm}((-1)^m R_n^{(m)})=m+2(n-1),\quad n\geq 2.\end{equation}
\begin{proposition}
For all $n\geq 2$, the completely monotonic degree of the function $R_n(x)$ with respect to $x>0$ satisfies: $$2(n-1)\leq {\rm deg}^x_{cm}(R_n)< 2n-1.$$
\end{proposition}
{\bf Proof.}  In [Qi and Mansour], it is proved that
$$-R'_{n+1}(x)=\int_0^\infty f_n(t)e^{-x t}dt,$$
where $$f_n(t)=(-1)^n(\frac1t-\frac12\coth\frac t2+\sum_{k=1}^{n+1}\frac{B_{2k}}{(2k)!}t^{2k-1}).$$
Moreover, $R_{n+1}$ is a completely monotonic of degree at lest $n+1$, hence, $\lim_{x\to\infty}R_{n+1}(x)=0$. Hence for all $x>0$,
We have, $$R_{n+1}(x)=\int_0^\infty g_n(t)e^{-xt}dt,$$
where $$g_n(t)=(-1)^n(\frac1{t^2}-\frac1{2t}\coth\frac t2+\sum_{k=1}^{n+1}\frac{B_{2k}}{(2k)!}t^{2k-2}).$$
Integrate by part yields $$x^{2n}R_{n+1}(x)=\int_0^\infty (g_n(t))^{(2n)}e^{-xt}dt,$$
\begin{equation}\label{i}g^{(2n)}_n(t)=(-1)^n((\frac1{t^2}-\frac1{2t}\coth \frac t2)^{(2n)}+\frac{B_{2n+2}}{(2n+2)(2n+1)}).\end{equation}
We use the Legendre integral formula. See for instance \cite{bi} (page 265) and \cite{et} (page 92).
 For all $t>0$,
$$\int_0^\infty \frac{\sin(x t)}{e^x-1}dx=\frac\pi2\coth (\pi t)-\frac1{2t}.$$
Integrate by part yields,
$$\int_0^\infty \frac{\sin(x t)}{e^x-1}dx=[\frac12\sin(xt)\log(1+e^{-2x}-2e^{-x})]_0^\infty-\frac t2\int_0^\infty\cos(xt)\log(1+e^{-2x}-2e^{-x})dx.$$
Hence,
$$\frac1{(2\pi)^2}\int_0^\infty\cos(\frac{xt}{2\pi})\log(1+e^{-2x}-2e^{-x})dx=\frac1{t^2}-\frac1{2t}\coth (\frac t2).$$
Applying the theorem of derivation under the integral sign, it follows that
$$g^{(2n)}_n(t)=\frac1{(2\pi)^{2n+2}}\int_0^\infty x^{2n}\cos(\frac{xt}{2\pi})\log(1+e^{-2x}-2e^{-x})dx+\frac{(-1)^nB_{2n+2}}{(2n+2)(2n+1)}.$$
Since, $$\frac1{t^2}-\frac1{2t}\coth(\frac  t2)=-\sum_{k=0}^\infty\frac{B_{2k+2}}{(2k+2)!}t^{2k},\qquad |t|<2\pi.$$
We deduce that,
 $$\frac1{(2\pi)^{2n+2}}\int_0^\infty x^{2n}\log(1+e^{-2x}-2e^{-x})dx=(-1)^{n+1}\frac{2^{2n}B_{2n+2}}{(2n+2)(2n+1)}.$$
 Thus,

$$g^{(2n)}_n(t)=\frac1{(2\pi)^{2n+2}}\int_0^\infty x^{2n}\Big(\cos(\frac{xt}{2\pi})-1\Big)\Big( \log(1+e^{-2x}-2e^{-x}))\Big)dx.$$
Let $\theta(x)= \log(1+e^{-2x}-2e^{-x})$. Then, $\theta'(x)=2(e^{-x}-e^{-2x})\geq 0$ and $\lim_{x\to\infty}\theta(x)=0$. then, $g^{(2n)}_n(t)>0$ for all $t>0$. Which implies that ${\rm deg}_{\rm cm}^t(R_{n+1})\geq 2n$.

Assume $t^\alpha R_n(t)$ is completely monotonic. Then $\alpha\leq -\frac{tR_n'(t)}{R_n(t)}$ for all $t>0$.
Since,
$$\frac{tR_n'(t)}{R_n(t)}=\frac{t(\psi(t)-\log t)-\frac1{2}+\sum_{k=1}^{n+1}\frac{B_{2k}}{2k}\frac1{t^{2k-1}}}{\log\Gamma(t)-(t-\frac12)\log(t)+t+\frac12\log(2\pi)-\sum_{k=1}^{n+1}\frac{B_{2k}}{2k(2k-1)}\frac1{t^{2k-1}}},$$
Hence, $$\lim_{t\to 0}\frac{-tR_n'(t)}{R_n(t)}=2n+1.$$
Furthermore,
$$x^{2n+1}R_{n+1}(x)=\int_0^\infty g^{(2n+1)}_n(t)e^{-xt}dt,$$
and
$$g^{(2n+1)}_n(t)=\frac{-1}{(2\pi)^{2n+3}}\int_0^\infty x^{2n+1}\sin(\frac{xt}{2\pi})\log(1+e^{-2x}-2e^{-x}))dx.$$
Hence, $x^{2n+1}R_{n+1}(x)$ is not completely monotonic. Then,  ${\rm deg}_{\rm cm}^t(R_{n+1})\in[2n,2n+1)$ for all $n$.
\begin{proposition} There is $m_0\in\Bbb N$ such that for $m\geq m_0$, the function $(-1)^mx^{m-1}R_0^{(m)}(x)$ is not completely monotonic.
\end{proposition}
{\bf Proof.} We have $$R_0(x)=\log\Gamma(x)-(x-\frac12)\log x+x-\frac12\log(2\pi).$$
then, for all $x>0$ and all $m\geq 2$, $$(-1)^mR_0^{(m)}(x)=(-1)^m\psi^{(m-1)}(x)-\frac{(m-1)!}{2x^m}-\frac{(m-2)!}{x^{m-1}}.$$
$$(-1)^mR_0^{(m)}(x)=\int_0^\infty(\frac {t^{m-1}}{1-e^{-t}}-\frac {t^{m-1}}2-t^{m-2})e^{-xt}dt,$$
then, for all $m\geq 2$, $$(-1)^mx^{m-1}R_0^{(m)}(x)=\int_0^\infty(\Big(\frac {t^{m-1}}{1-e^{-t}}\Big)^{(m-1)}-\frac{(m-1)!}2)e^{-xt}dt.$$
Set $g_m(t)=\Big(\frac {t^{m-1}}{1-e^{-t}}\Big)^{(m-1)}-\frac{(m-1)!}2$. 
Assume that $x^{m-1}\varphi_m$ is completely monotonic for all $m$, then, $g_m(t)\geq 0$ 
and $f_m(t)\geq \frac{(m-1)!}2$. Furthermore, it is proved by Alzer that for all $t>0$, $\displaystyle\lim_{m\to\infty}\frac1{(m-1)!}f_m(t/(m-1))=s(t)$, where,
$$s(t)=\frac12+\frac1\pi\sum_{k=1}^\infty\frac1k\sin(\frac t{2k\pi}).$$
Which gives, $s(t)\geq\frac12$ for all $t>0$. From Theorem 2.1
(see (\cite{alz1} p. 105)), we can derive that there is $a > 0$ such that $s(a) < 0$, Which gives a contradiction.
The Bernstein-Widder theorem \cite{sch} implies $x^{m-1}\varphi_m(x)$ is not completely monotonic on $(0,\infty)$ for all $m\in\Bbb N$. This completes the proof.
\begin{remark} By the proposition above, one deduces that $x^m(-1)^mR_0^{(m)}(x)$ is not completely monotonic for all $m\in\Bbb N$.
\end{remark}
One shows that, $m\geq 80$, $f'_m\Big(\sqrt{\frac{252}{(m+4)(m+3)}}\Big)<0$.
\begin{proposition} For $m\geq 3$, the completely monotonic degree of the function $(-1)^mR_0^{(m)}(x)$ with respect to $x>0$ is not less that $m-2$ and less than $m-1$,
$$m-2\leq{\rm deg}_{cm}^x((-1)^mR_0^{(m)}(x))<m-1$$
\end{proposition}
{\bf Proof.} As above, we have
$$x^{m-2}\varphi_m(x)=\int_0^\infty(\Big(\frac {t^{m-1}}{1-e^{-t}}\Big)^{(m-2)}-\frac{(m-1)!}2t-(m-2)!)e^{-xt}dt.$$
Moreover for $|x|<2\pi$,
$$\frac{x^{m-1}}{1-e^{-x}}=x^{m-2}+\frac{x^{m-1}}2+\sum_{k=2}^\infty\frac{B_k}{k!}x^{k+m-2},$$
hence, the $m-2$ derivative of the function  $\displaystyle\frac{x^{m-1}}{1-e^{-x}}$ at $0$ is equal to $(m-2)!$.
This implies, \begin{equation}\label{r8}\Big(\frac {t^{m-1}}{1-e^{-t}}\Big)^{(m-2)}-\frac{(m-1)!}2t-(m-2)!=\int_0^t(f_m(t)-\frac{(m-1)!}{2})dt.\end{equation}
To show that $x^{m-2}\varphi_m(x)$ is completely monotonic, it sufficient to prove that  for all $t>0$, $$\displaystyle\Big(\frac {t^{m-1}}{1-e^{-t}}\Big)^{(m-2)}-\frac{(m-1)!}2t-(m-2)!\geq 0.$$
We saw that $$f_{m}(t)=(m-1)!+(m-1)!\sum_{k=1}^\infty e^{-kt}L_{m-1}(kt),$$
and by the fact that, $|L_m(x)|\leq e^{\frac x2}$ for all $x>0$. It then follows that
$$|\sum_{k=1}^\infty e^{-kt}L_{m-1}(kt)|\leq \frac{e^{-\frac t2}}{1-e^{-\frac t2}},$$
hence, $$(m-1)!(1-\frac{e^{-\frac t2}}{1-e^{-\frac t2}})\leq f_m(t),$$
It is easy seeing that $1-\frac{e^{-\frac t2}}{1-e^{-\frac t2}}\geq\frac12$ if and only if $t\geq 2\log 3$. Then, for all $ t\geq 2\log 3\simeq 2.19$
$$ f_m(t)-\frac{(m-1)!}2\geq 0.$$
Moreover, H. Alzer et al. \cite{alz1} (p.113) showed that for all $t\in(-2\pi,2\pi)$,
$$f_m(t)=\int_0^\infty s(t u)u^{m-1}e^{-u}du,$$
where, $$s(u)=\frac12+\frac1\pi\sum_{k=1}^\infty\frac1k\sin(\frac u{2k\pi}).$$
Let $t\in[0,2\pi)$, then,
$$\int_0^{t}(f_m(t)-\frac{(m-1)!}2)dt=\int_0^{t}(\int_0^\infty (s(t u)-\frac12)u^{m-1}e^{-u}du)dt.$$
then,
$$\begin{aligned}\int_0^{t}(f_m(t)-\frac{(m-1)!}2)dt&=\int_0^\infty\Big(\int_0^{t}\frac1\pi\sum_{k=1}^\infty\frac1k\sin(\frac {tu}{2k\pi})\Big)dt)u^{m-1}e^{-u}du\\&=4\int_0^\infty\Big(\sum_{k=1}^\infty\sin^2(\frac {tu}{4k\pi})\Big)u^{m-2}e^{-u}du\geq 0.\end{aligned}$$
This complete the proof.
 \begin{proposition} For $m\geq 1$, the completely monotonic degree of the function $(-1)^mR_1^{(m)}(x)$ with respect to $x>0$ is not less that $m$ and less than $m+1$,
 $$m\leq{\rm deg}_{cm}^x((-1)^mR_1^{(m)}(x))<m+1$$
\end{proposition}
{\bf Proof.} We saw that
$$-R_1(x)=\log\Gamma(x)-(x-\frac12)\log(x)+x+\frac12\log(2\pi)-\frac1{12 x},$$
So, for all $m\geq 1$, $$(-1)^mR_1^{(m)}(x)=(-1)^{m+1}\psi^{(m-1)}(x)+\frac{(m-2)!}{x^{m-1}}+\frac{(m-1)!}{2x^m}+\frac{m!}{12x^{m+1}},$$
Hence, $$(-1)^mR_1^{(m)}(x)=\int_0^\infty\Big(\frac{t^m}{12}+\frac{t^{m-1}}{2}+t^{m-2}-\frac{t^{m-1}}{1-e^{-t}}\Big)e^{-xt}dt.$$
It follows that
$$(-1)^mx^mR_1^{(m)}(x)=\int_0^\infty(\frac{m!}{12}-f_m'(t))e^{-xt}dt,$$
where $\displaystyle f_m(t)=\Big(\frac{t^{m-1}}{1-e^{-t}}\Big)^{(m-1)}$. It is known that for $|t|<2\pi$,
$$f_m(t)=\int_0^\infty s(tu)u^{m-1}e^{-u}du,$$ then
$$f'_m(t)=\int_0^\infty s'(tu)u^{m}e^{-u}du.$$
Since, for all $x\in\Bbb R$, $\displaystyle s'(x)=\frac1{2\pi^2}\sum_{k=1}^\infty\frac1{k^2}\cos(\frac x{2k\pi})$. Then,
$|s'(x)|\leq \frac1{12}$. Thus,
\begin{equation}\label{k1}\frac{m!}{12}-f_m'(t)=\int_0^\infty(\frac1{12}-s'(tu))t^{m}e^{-t}dt\geq 0,\qquad{\rm for\;all}\;|t|<2\pi.\end{equation}
On the other hand,
$$f_m(t)=(m-1)!\sum_{k=0}^\infty e^{-kt}L_m(kt),$$
then, $$f_m'(t)=(m-1)!\sum_{k=0}^\infty ke^{-kt}(L'_m(kt)-L_m(kt)).$$
Using the relation, $tL_m'(t)=mL_m(t)-mL_{m-1}(t)$, then,
$$tf_m'(t)=(m-1)!\sum_{k=0}^\infty e^{-kt}(mL_m(kt)-mL_{m-1}(kt)-ktL_m(kt)).$$
By an inequality due to Szegö (see [14, p. 168]) we have $|L_m(t)|\leq e^{\frac t2}$ for $t\geq 0$, so that
we obtain for $t >0$,
$$tf_m'(t)\leq 2\,m!\frac{e^{-t/2}}{1-e^{-t/2}}-2t\,(m-1)!\frac{d}{dt}(\frac{e^{-t/2}}{1-e^{-t/2}}),$$
This yields the following inequality
\begin{equation}\label{p}tf_m'(t)\leq 2\,m!\frac{e^{-t/2}}{1-e^{-t/2}}+\,(m-1)!\frac{t}{4\sinh^2(t/4)},\end{equation}
which gives,
$$\frac{m!}{12}-f_m'(t)\geq (m-1)!(\frac m{12}-\frac{2me^{-t/2}}{t(1-e^{-t/2})}-\frac{1}{4\sinh^2(t/4)}),$$
Let $K(t,m)=\frac m{12}-\frac{2me^{-t/2}}{t(1-e^{-t/2})}-\frac{1}{4\sinh^2(t/4)}$.
It easy to see that the function $K(t,m)$ increases on the variable $m$ if and only if $\frac 1{12}-\frac{2e^{-t/2}}{t(1-e^{-t/2})}\geq 0$, which is true for $t\geq 4$.
Moreover, $$K'(t,1)=\frac{-2+e^{t/2} (2+t)}{(-1+e^{t/2})^2 t^2}+\frac18\frac{\coth(t/4)}{\sinh^2(t/4)}\geq 0.$$
 Then, for $t\geq 4$ and $m\geq 1$, $$K(t,m)\geq K(t,1).$$ Furthermore, for $t\geq 6$ we have $K(t,1)\geq K(6,1)>0.1$.
 Which implies that \begin{equation}\label{k2}K(t,m)\geq 0\qquad{\rm for\; all}\; t\geq 6,\;{\rm and}\;m\geq 1.\end{equation}
 By equations \eqref{k1} and \eqref{k2}, we get  for all $t>0$, $$\frac{m!}{12}-f_m'(t)\geq 0.$$
 Thus $(-1)^mx^mR_1^{(m)}(x)$ is completely monotonic.

 Let $m\geq 1$,
we have seen that $f'_m(t)\leq m!/12$ for all $t\geq 0$ and $f'_m(0)=m!/12$. Hence, $\displaystyle\lim_{t\to +\infty}(m!/2-f'_m(t))e^{-xt}=\displaystyle\lim_{t\to 0}(m!/2-f'_m(t))e^{-xt}=0$. Integrate by part yields
 $$(-1)^mx^{m+1}R_1^{(m)}(x)=-\int_0^\infty f_m''(t)e^{-xt}dt,$$
 If for all $m\geq 1$ and all $t>0$ $f''_m(t)\leq 0$. Then, by using equation \eqref{p}, we have,  $\lim_{t\to\infty} f'_m(t)=0$, and  $f'_m(t)\geq 0$ for all $t>0$. Therefore
 $$f_m(t)\geq f_m(0)=\frac{(m-1)!}{2}.$$
 Using the fact that $\displaystyle\lim_{m\to\infty}\frac1{(m-1)!}f_m(t/(m-1))=s(t)$ for all $t\in\Bbb R$. It follows that,
 $s(t)\geq 1/2$,  and  $H(t)\geq 0$ for all $t>0$.
Which contradicts the result of Alzer et al \cite{alz1}. Which states that $H(x_{j_k}) < -C(\log \log x_{k_j})^{1/2}$, $C>0$, for some positive sequence $x_{j_k}$ going to infinity as $k\to +\infty$.

The Bernstein-Widder theorem [36] implies $x^{m+1}(-1)^mR^{(m)}_1(x)$ is not completely monotonic on $(0,\infty)$ for all $m\in\Bbb N$. This completes the proof.

Address:  Institu pr\'eparatoire aux \'etudes d'ing\'enieurs de Tunis.\\
Campus Universitaire El-Manar, 2092 El Manar Tunis.\\
Email: bouali25@laposte.net
 \end{document}